\documentclass[11pt]{article}

\usepackage{amssymb,amsmath,amsfonts,amssymb}
\usepackage{graphics,graphicx,color,verbatim}

\def\@abssec#1{\vspace{.05in}\footnotesize \parindent .2in
{\bf #1. }\ignorespaces}

\graphicspath{{/EPSF/}{Figures/}}
\DeclareGraphicsExtensions{.eps}

\def \Rm {\mathbb R}

\newcommand{\dsum}{\displaystyle\sum}

\newcommand{\cout}[1]{}

 \renewcommand{\arraystretch}{1.5}

\title{Reconstructions for some coupled-physics inverse problems}

\author{Guillaume Bal \thanks{Department of Applied Physics and 
        Applied Mathematics, Columbia University, 
        New York NY, 10027; gb2030@columbia.edu} \and Gunther Uhlmann \thanks{Department of Mathematics, University of Washington, Seattle, WA, 98195 and University of California, Irvine, CA, 92697; gunther@math.washington.edu}}

\begin{document}
 
\maketitle


\begin{abstract}
   This letter announces and summarizes results obtained in \cite{BU-CPAM-12} and considers several natural extensions. The aforementioned paper proposes a procedure to reconstruct coefficients in a second-order, scalar, elliptic equation from knowledge of a sufficiently large number of its solutions. We present this derivation and extend it to show which parameters may or may not be reconstructed for several hybrid (also called coupled physics) imaging modalities including photo-acoustic tomography, thermo-acoustic tomography, transient elastography, and magnetic resonance elastography. Stability estimates are also proposed.
\end{abstract}
 

\renewcommand{\thefootnote}{\fnsymbol{footnote}}
\renewcommand{\thefootnote}{\arabic{footnote}}

\renewcommand{\arraystretch}{1.1}





\section{Introduction}
\label{sec:intro}

Consider a general second-order, linear elliptic equation of the form
\begin{equation}
\label{eq:elliptic}
\nabla\cdot a\nabla u_j + b\cdot\nabla u_j + cu_j =0,\quad x\in X,\qquad u_j=f_j \quad x\in \partial X,\qquad 1\leq j\leq J,
\end{equation}
for $X$ a smooth open domain in $\Rm^n$, with $n$ spatial dimension, and $(a,b,c)$ possibly complex-valued, symmetric second-order tensor, vector field, and scalar coefficient, respectively. We assume that $a$ is elliptic, the real part of $a$ is coercive and bounded, and $c$ is such that the above equation admits a unique solution. We also assume that $(a,b,c,\nabla\cdot a)$ are of class $C^{0,\alpha}$ for some $\alpha>0$. We construct $J$ solutions of the above equation for different boundary conditions. 

Several recent hybrid inverse problems aim to reconstruct the unknown coefficients $(a,b,c)$ from knowledge of internal functionals of the coefficients and of the elliptic solutions $(u_j)_{1\leq j\leq J}$. Concretely, we assume knowledge of the following functionals
\begin{equation}
\label{eq:functionals}
H_j(x) = d(x) u_j(x) ,\qquad x\in X,
\end{equation}
with $d(x)$ a scalar coefficient that is a priori also unknown. 

What may be reconstructed from $(a,b,c)$ from knowledge of $(H_j)_{1\leq j\leq J}$ in the setting $d\equiv1$ is analyzed in \cite{BU-CPAM-12}. We present the reconstruction procedure of the aforementioned reference in section \ref{sec:derivation}. Such a reconstruction is based on the availability of ratios of solutions $\frac{H_j}{H_k}=\frac{u_j}{u_k}$. This preliminary step is then used in section \ref{sec:reconstruction} to show that $(a,b,c,d)$ can be reconstructed up to explicit obstructions that take the form of gauge transformations. We also provide stability estimates for the reconstructions.

The inverse problems with internal functionals of the form \eqref{eq:functionals} are part of a larger class referred to as hybrid inverse problems or coupled-physics inverse problems. For recent results and reviews on this topic, we refer the reader to \cite{A-Sp-08,B-IO-12,CLB-SPIE-09,KK-EJAM-08,KLFA-MP-95,MZM-IP-10,OCPYL-UI-91,SU-IP-09,S-SP-2011,WW-W-07} and their references. 

\section{Reconstruction procedure}
\label{sec:derivation}

For the rest of the paper, we assume the existence of $u_1\not=0$ on $\bar X$. We refer to \cite{BU-CPAM-12} for conditions on $f_1$ that ensure such a property, either globally, in favorable cases, or at least locally. This allows us to define the known quantities
\begin{equation}
\label{eq:vj}
v_j=\dfrac{H_{j+1}}{H_1} = \dfrac{u_{j+1}}{u_1},\qquad 1\leq j\leq J-1.
\end{equation}
Using the notation $A:B={\rm Tr}(AB)$ for symmetric matrices $A$ and $B$, we find that 
\begin{equation}
\label{eq:forvj}
\alpha_\tau : \nabla ^{\otimes 2} v_j + \beta_\tau \cdot\nabla v_j =0 \quad x\in X,\qquad v_j=\frac{f_{j+1}}{f_1} \quad x\in \partial X,
\end{equation}
where for an arbitrary complex-valued non-vanishing function $\tau(x)$ on $X$, we have
\begin{equation}
\label{eq:alphabeta}
\alpha_\tau = \tau u_1^2 a,\qquad \beta_\tau = \tau u_1^2 b +\tau \nabla \cdot au_1^2. 
\end{equation}
Note that the equation \eqref{eq:forvj} is invariant by multiplication by a non-vanishing scalar coefficient so that $(\alpha_\tau,\beta_\tau)$ may at best be reconstructed up to a multiplicative scalar coefficient. The result in \cite{BU-CPAM-12} shows that this is the only obstruction to the reconstruction of $(\alpha_\tau,\beta_\tau)$.

More precisely, let us assume that $(\nabla v_1,\ldots,\nabla v_n)$ form a basis of $\Rm^n$ for all $x\in\bar X$. We distinguish the case $a$ scalar from the case $a$ a second-order tensor. When $a$ is scalar and $J=n+1$, then $(\alpha_\tau,\beta_\tau)$ are reconstructed up to the multiplicative scalar $\tau$. This is equivalent to saying that $a^{-1}b$ is uniquely reconstructed. Indeed, we have
\begin{displaymath}
\Delta v_j + \dfrac{\beta_\tau}{\alpha_\tau}\cdot \nabla v_j = \Delta v_j + \dfrac{b}{a}\cdot \nabla v_j =0
\end{displaymath}
so that, defining $H_{ij}=\nabla v_i\cdot\nabla v_j$ and $H^{ij}$ the coefficients of $H^{-1}$, we have
\begin{equation}
\label{eq:recb}
a^{-1}b = H^{ij} (a^{-1}b \cdot \nabla v_j) \nabla v_i = - H^{ij} \Delta v_j  \nabla v_i,
\end{equation}
where we have used the convention of summation over repeated indices and the fact that for any vector $F$, we have $F=H^{ij}F\cdot\nabla v_j \nabla v_i$.

When $a$ is tensor-valued, we need
\begin{equation}
\label{eq:numbermeas}
J = I_n := \dfrac12 n (n+3) = n+1+M_n ,\qquad M_n = \dfrac12n(n+1) -1.
\end{equation}
For $1\leq j\leq I_n-1$ and $1\leq m\leq M_n$, let us define the coefficients $\theta^m_j$ such that 
\begin{equation}
\label{eq:thetaM}
\dsum_{j=1}^{I_n-1} \theta^m_j \nabla v_j=0 \qquad \mbox{ and the symmetric matrices } \qquad 
M^m = \dsum_{j=1}^{I_n-1} \theta_j^m \nabla^{\otimes 2} v_j,
\end{equation}
such that $(M^m)_{1\leq m\leq M_n}$ form a free family of symmetric matrices. Sufficient conditions are presented in \cite{BU-CPAM-12} to guaranty that $(\nabla v_j)_{1\leq j\leq n}$ and $(M^m)_{1\leq m\leq M_n}$ are free families for the choice $\theta^m_j=-H^{jk}\nabla v_{m+n}\cdot\nabla v_k$ for $1\leq j\leq n$, $\theta^m_j=1$ for $j=n+m$ and $\theta^m_j=0$ otherwise. The above construction allows us to obtain the following constraints:
\begin{equation}
\label{eq:orthogalpha}
\alpha_\tau: M^m =0 ,\quad 1\leq m\leq M_n.
\end{equation}
This implies that $\alpha_\tau=M^0$, where $(M^0)^*$ is a matrix in the one-dimensional orthogonal complement to $(M^m)_{1\leq m\leq M_n}$ for the inner product for symmetric matrices $(A,B)={\rm Tr}(AB^*)$. Thus $\alpha_\tau$ is reconstructed up to a multiplicative scalar coefficient. From \eqref{eq:forvj}, we deduce that 
\begin{equation}
\label{eq:betatau}
\beta_\tau = -H^{ij} \alpha_\tau : \nabla^{\otimes 2} v_j \nabla v_i.
\end{equation}
Note that the above is nothing but \eqref{eq:recb}  when $a$ is a scalar coefficient.

This shows that $(\alpha_\tau,\beta_\tau)$ are uniquely reconstructed up to the multiplicative coefficient $\tau$. Note that additional information of the form $H_k=d u_k$ for $u_k$ solution of \eqref{eq:elliptic} with $u_k=f_k$ on $\partial X$ does not provide any new information. Indeed, $\frac{H_k}{H_1}$ is a solution of the elliptic equation $\eqref{eq:forvj}$ with known boundary condition $\frac{u_k}{u_1}$ on $\partial X$.

\section{Reconstruction of $(a,b,c,d)$ up to gauge transforms}
\label{sec:reconstruction}

\noindent{\bf Reconstruction up to gauge transforms.}
The above derivation shows that all that can be extracted from an arbitrary large number of functionals of the form $H_k=du_k$ is $(\alpha_\tau,\beta_\tau,H_1)$ augmented with the equation for $u_1$. Let us decompose $a=B^2\hat a$ for $\hat a$ a matrix with determinant equal to $1$. We assume here to simplify that such a decomposition is valid globally on $\bar X$ (which is obvious in the case where $a$ is real-valued and positive-definite). Since $\alpha_\tau=\tau u_1^2 B^2\hat a$ is known, we deduce that $\hat a$ is known. We compute
\begin{displaymath}
\hat a \alpha_\tau^{-1} (\beta_\tau-\nabla\cdot \alpha_\tau) = \dfrac{b}{B^2} - (\nabla \ln \tau)\cdot  \hat a.
\end{displaymath}
Moreover, defining $v=Bu_1=\frac{H_1B}d$, we find that 
\begin{displaymath}
 \dfrac{\Delta v}{v} = \dfrac{\nabla\cdot \hat a \nabla B}{B} + \dfrac{c}{B^2}.
\end{displaymath}
Thus, we obtain after elimination of $\tau$ and $u_1$ that knowledge of $(\alpha_\tau,\beta_\tau,H_1)$ and the equation for $u_1$ is equivalent to knowledge of
\begin{equation}\label{eq:coefs}
   \Big(\hat a \, , \,  \dfrac{b}{B^2} +2 \hat a \nabla \ln \dfrac B d\, , \, \dfrac{\Delta \frac{H_1B}{d}}{\frac{H_1B}{d}} = \dfrac{\nabla\cdot \hat a \nabla B}{B} + \dfrac{c}{B^2}\Big).
\end{equation}
No additional information may be extracted from functionals of the form $H_k=d u_k$ since knowledge of the above coefficients uniquely determines the functionals $H_k$.

The dimension of the unknown coefficients in \eqref{eq:coefs} is 
\begin{math}
  \frac{n(n+1)}2 -1 + n + 1 = I_n=\frac12 n(n+3),
\end{math}
which is the number of functionals used to reconstruct them. The dimension of $(a,b,c,d)$ is 
\begin{math}
  \frac{n(n+1)}2  + n + 1 +1 = I_n +2.
\end{math}
There are therefore two gauge parameters that remain undetermined. Moreover, $(a,b,c,d)$ are reconstructed up to any transformation that leave the coefficients in \eqref{eq:coefs} invariant.
\\[2mm]
\noindent{\bf Applications to medical imaging modalities.}
In the setting of Transient Elastography and Magnetic Resonance Elastography, we may assume that $d$ is known (and equal to $1$) and that $b=0$.  We thus obtain a (redundant) transport equation for $B$ (or equivalently for the gauge $\tau$) and then an explicit expression for $c$. Therefore, $(a,c)$ is uniquely reconstructed. More generally, when $\nabla\cdot (a^{-1}b)$ is known, we obtain an elliptic equation for $B$ or equivalently for $\tau$. Then $(a,a^{-1}b,c)$ is uniquely reconstructed. 

In the setting of quantitative photo-acoustic tomography (QPAT), we may assume that $b=0$ and that $d=\Gamma c$. We again obtain that $\frac{B}d=\frac{B}{\Gamma c}$ is known, and hence  $q=\frac{\Delta v}v$, is known. The reconstruction of $(B,c,\Gamma)$ is unique up to any transformation that leaves $(\frac{\Gamma c}B,\frac{\nabla\cdot \hat a \nabla B}{B} + \frac{c}{B^2})$ invariant. When $\Gamma$ is known, then $(B,c)$ are uniquely reconstructed \cite{BR-IP-11, BR-IP-12, BU-IP-10}.

A similar result may be obtained in the imaging modality called quantitative thermo-acoustic tomography (QTAT), where $d=\Gamma (\Im c) u_1^*$; see \cite{AGJN-QTAT-12,B-MED-12,BRUZ-IP-11} for a derivation of such a model for $H_j=\Gamma (\Im c) u_j u_1^*$. Assuming again that $b=0$, or more generally that $\nabla\cdot(a^{-1}b)$ is known so that $\tau$, or equivalently $\frac Bd$ is known, then $v=Bu_1$ is known. In this setting, we thus find that $(B,c,\Gamma)$ are reconstructed up to any transform that leaves  $(\frac{\Gamma \Im c}{B^2},\,\frac{\nabla\cdot \hat a \nabla B}{B} + \frac{c}{B^2})$ invariant. Note that when $a$ is real-valued, then $\Gamma$ is uniquely reconstructed and $(B,c)$ are reconstructed up to a transform that leaves $\frac{\nabla\cdot \hat a \nabla B}{B} + \frac{c}{B^2}$ invariant \cite{B-MED-12}.

Note that, more generally, one condition on the field $b$ is sufficient to uniquely reconstruct the gauge $\tau$ or equivalently $\frac Bd$. Indeed, we observe that the second known quantity in \eqref{eq:coefs} is equivalent to knowledge of $a^{-1}b+2\nabla\ln\frac Bd$. Thus, knowledge of one component of $a^{-1}b$, or of $\nabla\cdot a^{-1}b$, for instance, again provides an equation that allows us to uniquely reconstruct $\frac Bd$ and, hence, $a^{-1}b$. In such a setting, $q=\frac{\Delta v}v$ with $v=Bu_1=\frac {H_1B}{d}$ is known and $(B,c,d)$ can then be reconstructed up to any transform that leaves $(\frac{B}{d},\frac{\nabla\cdot \hat a \nabla B}{B} + \frac{c}{B^2})$ invariant.

\section{Sufficient conditions and stability estimates}
\label{stab}

\noindent{\bf Sufficient conditions.} The results of the preceding section exactly characterize which coefficients in $(a,b,c,d)$ can be reconstructed. Such reconstructions hinge on the solutions $(u_j)$ to be sufficiently independent. More precisely, we assume that $u_1\not=0$ on $\bar X$, $(\nabla v_j)_{1\leq j\leq n}$ is a basis of $\Rm^n$ at every point $x\in \bar X$, and that the matrices $M^m$ are linearly independent on $\bar X$. 

In some situations, for instance when complex geometric optics (CGO) solutions can be constructed, the above conditions are shown to hold for an open set of well-chosen boundary condition $(f_j)_{1\leq j\leq I_n}$ \cite{BU-CPAM-12}. However, in the general situation where $a$ is possibly complex-valued and anisotropic, such CGO solutions are not available. The linear independences mentioned above can be shown to hold locally on subdomains on $X$. More precisely, it is shown in \cite{BU-CPAM-12} that for a finite covering $\cup_{k=1}^K X_k$ of $X$, then for an open set of boundary conditions $(f_j)_{1\leq j\leq J=K\times I_n}$, we can construct a non-vanishing solution $u_{k,1}$ on $\bar X_k$, linearly independent gradients $(\nabla \frac{u_{k,j}}{u_{k,1}})_{2\leq j\leq n+1}$ and linearly independent matrices $M^{k,m}$ as constructed in \eqref{eq:thetaM}.
\\[2mm]
\noindent{\bf Stability estimates.} The procedure leading to the reconstruction of \eqref{eq:coefs} is explicit and allows one to estimate how errors in the functionals $(H_j)$ propagate into errors in the reconstructed coefficients. Let us assume that $d$ is known and smooth and that $b=0$ for concreteness. Similar results can be obtained in more general cases. We observe that the construction of the matrices $M^m$ involve taking two derivatives of the functionals $H_j$. The reconstruction of $\hat a$ therefore involves differentiating $(H_j)$ twice. 

When $b=0$, we observe that the reconstruction of $\nabla B$ or equivalently $\nabla \tau$ from the (second) vector field in \eqref{eq:coefs} also involves differentiating $(H_j)$ twice. Once $(\hat a,B)$ are known, then \eqref{eq:coefs} provides a formula for $c$. However, some simplifications occur. From \eqref{eq:alphabeta}, we observe that $\nabla\cdot a$ is reconstructed from differentiating $(H_j)$ twice (and not thrice). Then with $u_1$ known since $d$ is known, we reconstruct $c$ directly from \eqref{eq:elliptic} with again a loss of two derivatives. This yields the result
\begin{displaymath}
  \|(\hat a,c,\nabla\cdot a) -(\hat{\tilde a},\tilde c,\nabla\cdot\tilde a)\|_{C^{0,\alpha}} + \|B-\tilde B\|_{C^{1,\alpha}} \leq C \|(H_j-\tilde H_j)_{1\leq j\leq J}\|_{C^{2,\alpha}},
\end{displaymath}
for some positive constant $C$, where $\tilde H_j$ is constructed as $H_j$ in \eqref{eq:functionals} with the coefficients $(a,b,c)$ in \eqref{eq:elliptic} replaced by $(\tilde a,0,\tilde c)$. Similar stability estimates may be obtained in the more general case with $d$ and $b$ unknown; see \cite{BU-CPAM-12} for additional results.

\section*{Acknowledgment} 
GB was partially funded by grants NSF DMS-1108608 and DMS-0804696. GU was partially funded by the NSF and a Rothschild Distinguished Visiting Fellowship at the Newton Institute. 

{\footnotesize

}
\end{document}